\documentclass[12pt]{amsart}
\usepackage{amsmath,amssymb,amsthm}
\usepackage{amsfonts}
\usepackage[mathscr]{eucal}
\usepackage{graphicx}
\usepackage{amsmath,amssymb,amsthm}
\usepackage{amsfonts}
\usepackage[mathscr]{eucal}
\newtheorem{theorem}{Theorem}[section]
\newtheorem{lemma}[theorem]{Lemma}

\newtheorem{cor}[theorem]{Corollary}

\theoremstyle{definition}

\theoremstyle{remark}
\newtheorem{remark}[theorem]{Remark}

\numberwithin{equation}{section}

\newcommand{\rd}{{\mathbb R^d}}

\newcommand{\rr}{{\mathbb R}}
\newcommand{\alert}[1]{\fbox{#1}}

\def\a{\alpha}

\def\E{{\mathbb E}}

\begin{document}

\title{\bf Space-time duality for fractional diffusion}
\author{Boris Baeumer}
\address{Boris Baeumer, Department of Mathematics and Statistics, University of Otago, PO. Box 56,
Dunedin, NZ} \email{bbaeumer@maths.otago.ac.nz}

\author{Mark M. Meerschaert}
\address{Mark M. Meerschaert, Department of Probability and Statistics,
Michigan State University, East Lansing, MI 48823}
\email{mcubed@stt.msu.edu}
\urladdr{http://www.stt.msu.edu/$\sim$mcubed/}
\thanks{Research of M. M. Meerschaert was partially
supported by NSF grant DMS-0706440.}

\author{Erkan Nane}
\address{Erkan Nane, 221 Parker Hall, Department of Mathematics and Statistics,
Auburn University, Auburn, Al 36849}
\email{nane@auburn.edu}
\urladdr{http://www.auburn.edu/$\sim$ezn0001}

\begin{abstract}
Zolotarev proved a duality result that relates stable densities with different indices.  In this paper, we show how Zolotarev duality leads to some interesting results on fractional diffusion.  Fractional diffusion equations employ fractional derivatives in place of the usual integer order derivatives.  They govern scaling limits of random walk models, with power law jumps leading to fractional derivatives in space, and power law waiting times between the jumps leading to fractional derivatives in time.  The limit process is a stable L\'evy motion that models the jumps, subordinated to an inverse stable process that models the waiting times.  Using duality, we relate the density of a spectrally negative stable process with index $1<\alpha<2$ to the density of the hitting time of a stable subordinator with index $1/\alpha$, and thereby unify some recent results in the literature.  These results also provide a concrete interpretation of Zolotarev duality in terms of the fractional diffusion model.
\end{abstract}


\maketitle


\section{Introduction}
A classical result of Zolotarev \cite{zolotarev-61} (see also Lukacs \cite[Theorem 3.3]{lukacs}) equates stable densities with different indices.  The proof of this result is purely analytical.  In this paper, we apply Zolotarev duality to prove some interesting results on fractional diffusion.  Fractional derivatives are natural extensions of their integer order analogues \cite{MillerRoss,Samko}.  Partial differential equations of fractional order are important in applications to physics \cite{MetzlerKlafter}, finance \cite{scalas1}, and hydrology \cite{FADEreview}.  In some cases, the solutions of the fractional equations govern the probability densities of certain heavy tailed stochastic processes \cite{Zsolution,limitCTRW,OB,OB1}. This connection, a generalization of the link between Brownian motion and the diffusion equation \cite{Einstein1906}, is very useful in both theoretical and applied work \cite{Bhattacharya,Zhang-PRE2008}.  Perhaps the simplest version of the fractional diffusion equation is $\partial^\gamma u/\partial t^\gamma = \partial^2 u /\partial x^2$ in one dimension, where the usual first derivative in time is replaced by the Caputo fractional derivative of order $0<\gamma<1$.  Meerschaert et al. \cite{Zsolution,limitCTRW} shows that the point source solution $u(x,t)$ to this equation gives the density of the stochastic process $B(E_t)$, where $B(x)$ is a Brownian motion and $E_t$ is the inverse or hitting time of a stable subordinator with index $\gamma$.  Orsingher and Beghin \cite{OB,OB1} show that the same solution can be written in terms of the normal density of $B(x)$ subordinated to a stable density with index $\alpha=1/\gamma$.  In this paper, we reconcile these two results using Zolotarev duality.  As a consequence, we reveal a concrete interpretation of the duality in terms of stable processes and their inverses.

\section{Duality}\label{sec2}
Stable laws are important because they represent the most general distributional limit for centered and normalized sums of independent and identically distributed random variables \cite{Feller1971,RVbook}.  Since most stable densities cannot be written in closed form, it is common to use characteristic functions (Fourier transforms).  A stable density $p(x)$ has characteristic function $\hat p(\lambda)=\int_{-\infty}^\infty e^{i\lambda x} p(x)\,dx$ and $p(x)=\frac1 {2\pi} \int_{-\infty}^\infty e^{-i\lambda x} \hat p(\lambda)\,d\lambda$.  Several different parameterizations of the family of stable densities are commonly used in the literature.   One commonly used representation of the centered stable density is \cite[Theorem 2.3]{lukacs}
\begin{equation}\label{common-char-func}
p_{\alpha}(x;\theta, c)=\frac1 {2\pi} \int_{-\infty}^\infty e^{-i\lambda x}\exp\bigg\{ -c |\lambda|^{\alpha}\bigg[1+i\theta\frac\lambda {|\lambda|}\tan (\pi\alpha/2)\bigg]\bigg\} d\lambda, \ \ \alpha\neq 1
\end{equation}
where $c>0$, $|\theta|\leq 1$, $0<\alpha\leq 2$, and  $\alpha\neq 1$.

A second parametrization is \cite{Feller1971,lukacs}:
\begin{equation}\label{char-func}
p_{\alpha}(x;\eta, b)=\frac1 {2\pi} \int_{-\infty}^\infty e^{-i\lambda x}\exp\bigg\{ -b |\lambda|^{\alpha}e^{-\frac{i\pi\eta}{2}\frac\lambda {|\lambda|}}\bigg\} d\lambda, \ \ \alpha\neq 1
\end{equation}
were $b>0$ and $\eta$ is real.
The connection between \eqref{common-char-func} and \eqref{char-func} is as follows:
$$c=b \cos\big( \frac{\pi\eta}{2}\big)\quad{\rm and}\quad\theta=-\cot \big(\frac{\pi \alpha}{2}\big)\tan \big(\frac{\pi \eta}{2}\big).$$
From the relation $|\theta|\leq 1$ one can conclude that $|\eta|\leq \alpha$ if $0<\alpha<1$ but $|\eta|\leq 2-\alpha$ if $1<\alpha\leq 2$.

Finally, the commonly used parametrization of Samorodnitsky and Taqqu \cite[Definition 1.1.6]{ST94} is
\begin{equation}\label{STparam}
p_{\alpha}(x;\beta, \sigma)=\frac1 {2\pi} \int_{-\infty}^\infty e^{-i\lambda x}\exp\bigg\{ -\sigma^\alpha |\lambda|^{\alpha}\bigg[1-i\beta\frac\lambda {|\lambda|}\tan (\pi\alpha/2)\bigg]\bigg\} d\lambda, \ \ \alpha\neq 1
\end{equation}
which is very similar to \eqref{common-char-func} except for a sign change.  The scale parameter $\sigma>0$ satisfies $c=\sigma^\alpha$ and the parameter $\beta=-\theta$ is called the skewness.  The change of sign between \eqref{common-char-func} and \eqref{STparam} has led to some confusion in the literature.

A duality result for stable densities was proven by Zolotarev \cite{zolotarev-61} (see also Lukacs \cite{lukacs}) using the parametrization \eqref{char-func}.  The proof follows directly from the series representation for stable densities \cite[Theorem 3.1]{lukacs}.

\begin{theorem}\label{duality-thm}
If $1<\alpha\leq 2$ then for all $u>0$ we have
\begin{equation}\label{duality}
p_\a(u;\eta,1)= u^{-(1+\alpha)}p_{\alpha^*}(u^{-\alpha}; \eta^*, 1)
\end{equation}
where $\alpha^*=1/\alpha$ and $\eta^*=\frac{\eta-1}{\alpha}+1$.
\end{theorem}

We will be interested in the case where the $\alpha$-stable density on the left-hand side of \eqref{duality} is totally negatively skewed ($\beta=-1$) which corresponds to $\theta=+1$ or equivalently $\eta=2-\alpha$.  In that case, $\eta^*=\alpha^*$ and so the stable density on the right-hand side has $\theta^*=-1$, or in other words, its skewness is $\beta^*=+1$ (totally positively skewed).  Then the right-hand side of \eqref{duality} involves the density of a stable subordinator whose index $\gamma=\alpha^*=1/\alpha$ satisfies $1/2 \leq\gamma<1$.  After substituting back into \eqref{char-func}, a little algebra shows that the characteristic function of the subordinator density is $\exp(-(-i\lambda)^\gamma)$.

\section{Fractional diffusion}\label{sec3}
Fractional diffusion equations are the governing equations of certain stochastic processes that occur as scaling limits of continuous time random walks \cite{Zsolution,limitCTRW}.  A continuous time random walk (CTRW) is simply a random walk in which the IID jumps $(Y_n)$ are separated by IID random waiting times $(J_n)$.  The CTRW was developed as a model in statistical physics \cite{MontrollWeiss,ScherLax}.  A random particle jump $Y_n$ follows a random waiting time $J_n>0$.  The random walk $S(n)=Y_1+\cdots+Y_n$ gives the particle location after $n$ jumps.  Another random walk $T(n)=J_1+\cdots+J_n$ gives the time of the $n$th jump.  The number of jumps by time $t>0$ is given by the renewal process $N_t=\max\{n\geq 0:T(n)\leq t\}$.  The location of a particle at time $t>0$ is $S(N_t)$, a random walk subordinated to a renewal process.  The long-time behavior of particles is described by a limit theorem \cite{limitCTRW}.  Suppose that $P(J_n>t)=t^{-\gamma}L_1(t)$ for some $0<\gamma<1$ where $L_1$ is slowly varying.  Then $T_n$ belongs to the domain of attraction of some stable law with index $\gamma$.  To simplify the exposition, consider the special case where $L_1(t)\to C>0$ is asymptotically constant as $t\to\infty$.  Then $n^{-1/\gamma}T_n\Rightarrow D$ as $n\to\infty$ where $D$ is a stable random variable with index $\gamma$.  Suppose also that $Y_n$ belongs to the domain of attraction of some stable random variable $A$ with index $0<\delta\leq 2$ and that $n^{-1/\delta}S_n\Rightarrow A$.  Restricting to the uncoupled case where $Y_n,J_n$ are independent (see \cite{coupleCTRW} for the coupled case) we can extend to process convergence $c^{-1/\gamma}T_{[ct]}\Rightarrow D(t)$ and $c^{-1/\delta}S_{[ct]}\Rightarrow A(t)$ where $A(t),D(t)$ are independent L\'evy stable processes with $A(1)=A,D(1)=D$ in distribution, and the convergence is in terms of all finite dimensional distributions (or the appropriate Skorokhod topology, see \cite{limitCTRW}).  Define the inverse or hitting time process
\begin{equation}\label{Edef}
E_t=\inf\{x>0:D(x)>t\}
\end{equation}
and note that $\{D(x)>t\}=\{E_t<x\}$.  Then $c^{-\gamma}N_{[ct]}\Rightarrow E_t$ \cite[Theorem 3.2]{limitCTRW} and a continuous mapping argument \cite[Theorem 4.2]{limitCTRW} shows that $c^{-\gamma/\delta} S(N_{[ct]})\Rightarrow A(E_t)$.

In applications to differential equations, the following parametrization of stable densities is useful \cite[Remark 11.1.13]{RVbook}:
\begin{equation}\label{FPDE-char-func}
p_{\alpha}(x;q,a)=\frac1 {2\pi} \int_{-\infty}^\infty e^{-i\lambda x}\exp\bigg\{ qa (i\lambda)^\alpha + (1-q)a(-i\lambda)^\alpha\bigg\} d\lambda, \ \ \alpha\neq 1
\end{equation}
which is related to \eqref{STparam} by $\beta=1-2q$ and $\sigma^\alpha=-a\cos(\pi\alpha/2)$.  Note that $0\leq q \leq 1$, and $a>0$ for $1<\alpha<2$, while $a<0$ for $0<\alpha<1$.  If $A=A(1)$ has density $p_{\delta}(x;q,a)$, then $A(t)$ has density $p_{\delta}(x;q,at)$.  Define $d^\delta p(x)/d x^\delta$ to be the function with Fourier transform $(-i\lambda)^\delta \hat p(\lambda)$, which extends the familiar formula for integer order derivatives.  Similarly, let $d^\delta p(x)/d (-x)^\delta$ to be the function with Fourier transform $(i\lambda)^\delta \hat p(\lambda)$.  These are called the (positive and negative) Riemann-Liouville fractional derivatives.  Since
\[\hat p(\lambda,t)=\exp\bigg\{ qat (i\lambda)^\delta + (1-q)at(-i\lambda)^\delta\bigg\}\]
is evidently the solution to the ordinary differential equation
\[\frac d{dt}\hat p(\lambda,t)=\bigg\{ qa (i\lambda)^\delta + (1-q)a(-i\lambda)^\delta\bigg\}\hat p(\lambda,t)\]
with the point source initial condition $\hat p(\lambda,0)\equiv 1$, we can invert the Fourier transform to see that
\begin{equation}\label{fADE-1L}
\frac{\partial p(x,t)}{\partial t}=qa\frac{\partial^\delta p(x,t)}{\partial
(-x)^\delta}+(1-q)a\frac{\partial^\delta p(x,t)}{\partial x^\delta}
\end{equation}
the space-fractional diffusion equation \cite{Chaves}.  We also call \eqref{fADE-1L} the governing equation of the process $A(t)$.  Note that $A(t)$ satisfies a tail condition $P(|A(t)|>x)\sim Ct x^{-\delta}$ where $C=-a/\Gamma(1-\delta)$ for $0<\delta<1$ and $C=a(\delta-1)/\Gamma(2-\delta)$ for $1<\delta<2$ \cite[Proposition 1.2.15]{ST94}.  Thus the order of the fractional derivative equals the tail index of the stable law. In addition, recall that $n^{-1/\delta}S_n\Rightarrow A$ which requires $P(|Y_n|>x)\sim C x^{-\delta}$ for $x$ large, so the order of the fractional derivative also reflects the tail behavior of particle jumps.  Finally, note that $P(Y_n<-x)/P(|Y_n|>x)\to q$ as $x\to\infty$, so that the positive and negative fractional derivatives in the governing equation reflect the positive and negative tails of particle jumps.

The governing equation for the density $f(x,t)$ of the time process $D(t)$ is similar:
\begin{equation}\label{Dgov}
\frac{\partial f(x,t)}{\partial t}=-b\frac{\partial^\gamma p(x,t)}{\partial x^\gamma}
\end{equation}
where $b>0$, and only the positive fractional derivative appears, since $J_n>0$.  Note that $f(x,t)=p_\gamma(x;\gamma,bt)$ in the parametrization \eqref{char-func}.  Write $g_\gamma(x)=f(x,1)$ and note that $f(x,t)=t^{-1/\gamma}g_\gamma(t^{-1/\gamma}x)$ (self-similarity).  The inverse process $E_t=(D/t)^{-\gamma}$ in distribution, and it follows that $x=E_t$ has density \cite[Corollary 3.1]{limitCTRW}:
\begin{equation}\label{Edens}
h(x,t)=\frac{t}{\gamma}x^{-1-1/\gamma}g_\gamma(tx^{-1/\gamma})
\end{equation}
on $x>0$ for all $t>0$.  This density solves \cite[Theorem 4.1]{triCTRW}:
\begin{equation}\label{ICP}
\frac{\partial h(x,t)}{\partial x}=-b\frac{\partial^\gamma h(x,t)}{\partial t^\gamma} +b\delta(x)\frac{t^{-\gamma}}{\Gamma(1-\gamma)}
\end{equation}
where we note that the roles of space and time are reversed for this inverse process density.  The Laplace transform $\tilde h(x,s)=\int_0^\infty e^{-s t}h(x,t)\,dt$ exists for all $x>0$ \cite{triCTRW}, and we may also consider ${\partial^\gamma h(x,t)}/{\partial t^\gamma}$ as the inverse Laplace transform of $s^\gamma \tilde h(x,s)$.  Now a simple conditioning argument shows that the CTRW scaling limit $A(E_t)$ has density
\begin{equation}\label{h1}
m(x,t)=\int_0^\infty p(x,u)h(u,t)\,du
\end{equation}
which randomizes the density of $x=A(u)$ according to the hitting time process $u=E_t$.  The overall governing equation is
\begin{equation}\label{gCp}
b\frac{\partial^\gamma m(x,t)}{\partial t^\gamma} =qa\frac{\partial^\delta m(x,t)}{\partial
(-x)^\delta}+(1-q)a\frac{\partial^\delta m(x,t)}{\partial x^\delta}+b\delta(x)\frac{t^{-\gamma}}{\Gamma(1-\gamma)}
\end{equation}
using the Riemann-Liouville fractional derivatives on both sides.  The Caputo fractional derivative $(\partial/\partial t)^\gamma F(t)$, defined for $0<\gamma\leq 1$ as the inverse Laplace transform of $s^\gamma \tilde F(s)-s^{\gamma-1}F(0)$, extends the usual integer order formula, and is useful in differential equations since it includes the initial value.  Then we can write \eqref{gCp} more compactly in the form
\begin{equation}\label{gCpC}
b\left(\frac{\partial }{\partial t}\right)^\gamma m(x,t)=qa\frac{\partial^\delta m(x,t)}{\partial
(-x)^\delta}+(1-q)a\frac{\partial^\delta m(x,t)}{\partial x^\delta}
\end{equation}
since $s^{\gamma-1}$ is the Laplace transform of ${t^{-\gamma}}/{\Gamma(1-\gamma)}$.  Similarly, we can rewrite \eqref{ICP} in the form
\begin{equation}\label{ICPC}
b\left(\frac{\partial }{\partial t}\right)^\gamma h(x,t)=-\frac{\partial h(x,t)}{\partial x}
\end{equation}
which can be considered a degenerate case of \eqref{gCpC} with $x=A(u)=u$.

\section{Space-time duality}
Here we apply Zolotarev duality from Section \ref{sec2} to the space-time fractional diffusion equation from Section \ref{sec3}.  The following result uses the parametrization \eqref{char-func}.

\begin{theorem}\label{YDlink}
Let $x=D(t)$ be a stable subordinator with density $p_{\gamma}(x;\gamma, bt)$  for some $1/2\leq\gamma<1$.  Let $E_t$ be the hitting time \eqref{Edef} with density \eqref{Edens}, where $g_\gamma(x)=p_{\gamma}(x;\gamma, b)$ is the density of $D(1)$, and let $Y(t)$ denote a stable L\'evy motion with density $p_{\alpha}(x;2-\alpha, b^{-\alpha}t)$ where $\alpha=1/\gamma$. Then
\begin{itemize}
\item[(i)] $P(Y(t)>0)=1/\alpha \ \ $ for all $t>0$.
\item[(ii)] $E_t$ is identically distributed with $Y(t)| Y(t)>0$ for each $t>0$.
\end{itemize}
\end{theorem}

\begin{proof}
Using the self-similarity property $p_\alpha(u;\eta,b)=b^{-1/\alpha}p_\alpha(b^{-1/\alpha}u;\eta,1)$ for stable densities, the density of $Y(t)$ for $t>0$ is
\[P(x,t)=bt^{-1/\alpha}p_{\alpha}(bt^{-1/\alpha}x;\eta, 1) .\]
Apply \eqref{duality} with $u=bt^{-1/\alpha}x$ and $\eta=2-\alpha$ to see that
\[P(x,t)=bt^{-1/\alpha} u^{-1-\alpha}p_{\alpha^*}(u^{-\alpha}; \eta^*, 1)\]
for $x>0$, where $\alpha^*=1/\alpha=\gamma$ and $\eta^*=\alpha^{-1}(\eta-1)+1=\gamma$.  Simplify to get
\[P(x,t)=tx^{-1-1/\gamma}b^{-1/\gamma}p_{\gamma}(b^{-1/\gamma}tx^{-1/\gamma}; \gamma, 1)\]
and recall that $g_\gamma(x)=p_{\gamma}(x; \gamma, b)=b^{-1/\gamma}p_{\gamma}(b^{-1/\gamma}x; \gamma, 1)$ is the density of $D(1)$.  Then
\[P(x,t)=tx^{-1-1/\gamma}g_{\gamma}(tx^{-1/\gamma})\]
for $x>0$.  Compare with \eqref{Edens} to conclude that
\begin{equation}\label{EtoY}
P(x,t)=\gamma h(x,t)
 \end{equation}
for $t>0$ and $x>0$.  Since $\int_0^\infty h(x,t)\,dx=1$ for all $t>0$ we have $P(Y(t)>0)=\gamma=1/\alpha$ for all $t>0$.  Thus, the density of $E_t$ equals the conditional density of $Y(t)$ given $Y(t)>0$ for all $t>0$.
\end{proof}

\begin{remark}\label{BinghamRemark}
As the scaling limit of a random walk with positive jumps, the stable subordinator $D(t)$ is totally positively skewed with $\beta=1$ in \eqref{STparam}.  The process $Y(t)$ has skewness $\beta=-1$ so it is the scaling limit of a random walk with only negative jumps.  This is also called a spectrally negative process, since the L\'evy measure assigns no mass to the positive real line.  Bingham \cite{bingham} points out that the hitting time $D(t)=\inf\{u:Y(u)>t\}$ is a stable subordinator with index $1/\alpha$ and that there is a version of $Y(t)$ for which $Y(D(t))\equiv t$.  Hence $D(t)$ is the process inverse of $Y(t)$, but $E_t$ is the process inverse of $D(t)$.  Thus the inverse of the inverse of $Y(t)$ is the process $E_t$ whose one dimensional distributions are the same as those of $Y(t)|Y(t)>0$.  Note that $D(E_t)>0$ almost surely since $D(t)$ is a pure jump process \cite{bertoin}.  If $\alpha=2$ then $Y(t)$ is a Brownian motion, and the skewness is irrelevant.
\end{remark}

\begin{remark}
Using the series representation for stable densities \cite[Theorem 3.1]{lukacs} it is not hard to prove Theorem \ref{YDlink} directly.  For convenience we consider $b=t=1$; the remaining cases follow easily by self-similarity.  Then the density of $Y(1)$ for $x>0$ is
\begin{equation}\label{alpha-big}
p_{\alpha}(x;\eta, 1)=\frac{1}{\pi}\sum_{k=1}^\infty (-1)^{k+1}\frac{\Gamma (1+k/\alpha)}{k!}x^{k-1}\sin (\frac{\pi k(\eta+\alpha)}{2\alpha} )
\end{equation}
where $\eta=2-\alpha$.  The density of stable subordinator $D$ is
\begin{equation}
g_\gamma(x)=\frac{1}{\pi}\sum_{k=1}^\infty (-1)^{k+1}\frac{\Gamma (\gamma k +1)}{k!}x^{-\gamma k-1}\sin (\pi k \gamma)
\end{equation}
Then the density of $E_1$ is
\begin{equation}
\begin{split}
x^{-1-\alpha}g_\beta (x^{-\alpha})
&=\alpha \frac{x^{-1-\alpha}}{\pi}\sum_{k=1}^\infty (-1)^{k+1}\frac{\Gamma (1+k/\alpha)}{k!}x^{(-k/\alpha-1)(-\alpha)}\sin (\frac{\pi k}\alpha )\\
&=\alpha \frac{1}{\pi}\sum_{k=1}^\infty (-1)^{k+1}\frac{\Gamma (1+k/\alpha)}{k!}x^{k-1}\sin (\frac{\pi k}\alpha )\\
&=\alpha p_{\alpha}(x; (2-\alpha), 1) .
\end{split}
\end{equation}
\end{remark}

Theorem \ref{YDlink} has some immediate consequences for space-time diffusion equations.  The discussion in Section \ref{sec3} explains how space-fractional derivatives model heavy tailed power law particle jumps, and time-fractional derivatives model power law waiting times.  The duality Theorem \ref{YDlink} connects heavy tailed jumps with fractional time derivatives, and conversely, it relates heavy tailed waiting times to fractional derivatives in space.  The Caputo fractional derivative
\[\left(\frac{\partial}{\partial t}\right)^\gamma h(x,t)=\frac 1{\Gamma(1-\gamma)} \int_0^t \frac{\frac {d}{dr}h(x,r)}{(t-r)^{\gamma}}\,dr\]
was used in the governing equation \eqref{ICPC}.  Then the next result follows immediately from \eqref{EtoY}.

\begin{cor}
Let $Y(t)$ denote a stable L\'evy motion with index $1<\alpha\leq 2$ and density $P(x,t)=p_{\alpha}(x;2-\alpha, b^{-\alpha}t)$ in the parametrization \eqref{char-func}.  Then
\begin{equation}\label{Pgovtime}
b\left(\frac{\partial }{\partial t}\right)^\gamma P(x,t)=-\frac{\partial P(x,t)}{\partial x}
\end{equation}
holds for all $t>0$ and $x>0$.
\end{cor}

The process $Y(t)$ is totally negatively skewed, so its density $P(x,t)$ solves a space-fractional equation similar to \eqref{fADE-1L} with $q=1$, using a negative Riemann-Liouville fractional derivative
\begin{equation}\label{RLnegdef}
\frac{\partial^\alpha P(x,t)}{\partial (-x)^\alpha}=\frac 1{\Gamma(2-\alpha)} \frac {d^2}{dx^2}\int_x^\infty \frac{P(y,t)}{(y-x)^{\alpha-1}}\,dy
\end{equation}
with $1<\alpha<2$.  In general, the $\alpha$ order negative Riemann-Liouville fractional derivative is defined as the $n$th derivative of a fractional integral of order $n-\alpha$, where $n-1<\alpha<n$ \cite{Samko}.

\begin{cor}\label{h-xgov}
Let $D(1)$ be a stable subordinator with density $g_\gamma(x)=p_\gamma(x;\gamma,b)$ in the parametrization \eqref{char-func}.  Let $E_t$ denote the hitting time process defined by \eqref{Edef}.  Then the density $h(x,t)$ of $E_t$ solves
\begin{equation}\label{Egovspace}
\frac{\partial h(x,t)}{\partial t}=b^{-\alpha}\frac{\partial^\alpha h(x,t)}{\partial (-x)^\alpha}
\end{equation}
for all $t>0$ and $x>0$.
\end{cor}

\begin{proof}
A comparison with \eqref{fADE-1L} shows that the density $P(x,t)$ of $Y(t)$ solves the space-fractional diffusion equation
\begin{equation}\label{Ygovspace}
\frac{\partial P(x,t)}{\partial t}=b^{-\alpha}\frac{\partial^\alpha P(x,t)}{\partial (-x)^\alpha}
\end{equation}
where we note that $c=\sigma^\alpha=-a\cos(\pi\alpha/2)=b^{-\alpha}\cos(\pi(2-\alpha)/2)$.  Then \eqref{Egovspace} follows using \eqref{EtoY} and the fact that the negative fractional derivative in \eqref{Ygovspace} depends only on $P(y,t)$ for $y>x$.
\end{proof}

\begin{remark}
The characteristic function (Fourier transform) of $P(x,t)$ is $\hat P(\lambda,t)=\exp(tb^{-\alpha}(i\lambda)^\alpha)$.  It is common to analyze partial differential equations like \eqref{Ygovspace} using transforms.  Note, however, that the Fourier transform of $P(x,t)I(x>0)$ is not equal to $\hat P(\lambda,t)$ since $P(x,t)$ is supported on the entire real line.  Since we restrict to the positive reals, it is convenient to use Laplace transforms.  Bingham \cite{Bingham} and Bondesson, Kristiansen, and Steutel \cite{BKS} show that
\[\int_0^\infty e^{-zx} h(x,t)\,dx=E_\gamma(-zb^{-1}t^\gamma)\]
in terms of the Mittag-Leffler function
$$E_\beta(z)=\sum_{k=0}^{\infty}\frac{z^k}{\Gamma (1+\beta k)}. $$
Then it follows from \eqref{EtoY} that
\[\int_{-\infty}^\infty e^{-z x} P(x,t)I(x>0)\,dx=\gamma E_\gamma(-z b^{-1} t^\gamma) \]
which shows that the conditional distribution of $Y(t)|Y(t)>0$ is Mittag-Leffler.
\end{remark}

For an $\rd$-valued Markov process $X(t)$, the family of linear operators
$T(t)r(x)=\E_x[r(X(t))]=E[r(X(t))|X(0)=x]$ forms a bounded
continuous semigroup
on the Banach space $L^1(\rd)$, and the generator $L_x r(x)=\lim_{h\downarrow 0}
  h^{-1}(T(h)r(x)-r(x))$ is defined on a dense subset of that space \cite{ABHN,HiPh}.  Then
$p(x,t)=T(t)r(x)$ solves the abstract Cauchy problem
\begin{equation}\label{semigroup-0}
\frac{\partial}{\partial t}p(x,t) =L_x p(x,t); \quad p(x,0) = r(x)
\end{equation}
for $t>0$ and $x\in \rd$.  Nigmatullin \cite{nigmatullin} considered an abstract time-fractional Cauchy problem
\begin{equation}\label{frac-derivative-0}
\left(\frac{\partial }{\partial t}\right)^\gamma m(x,t) = {L_x}m(x,t); \quad m(0,x) = r(x)
\end{equation}
which reduces to \eqref{gCpC} in the special case where $X(t)$ is a stable L\'evy process started at $x=0$.  Zaslavsky \cite{zaslavsky} used \eqref{frac-derivative-0} to model Hamiltonian chaos.  Baeumer and Meerschaert \cite{fracCauchy} and Meerschaert and Scheffler \cite{limitCTRW} showed that the solution to \eqref{frac-derivative-0} can be written in the form \eqref{h1} where $h(x,t)$ is the density of $E_t$, the hitting time \eqref{Edef} of a standard stable subordinator $D(t)$, and $g_{\gamma}(x)=p_\gamma(x;\gamma,1)$ is the density of $D(1)$.  It follows easily that $m(x,t)=\E_x[r(X({E(t)}))]$.  Then the next results follows immediately using \eqref{EtoY}.

\begin{lemma}\label{Ysolve}
Let $Y(t)$ denote a totally negatively skewed stable L\'evy motion with index $1<\alpha\leq 2$ and density $p_{\alpha}(x;2-\alpha, t)$ in the parametrization \eqref{char-func}.  Then the abstract fractional Cauchy problem \eqref{frac-derivative-0} with $\gamma=1/\alpha$ can be solved by taking
\[m(x,t)=\E_x[r(X(Y(t)))|Y(t)>0].\]
\end{lemma}

Next we come to the problem that motivated this paper.  For $\frac12\leq \gamma <1$, Orsingher and Beghin \cite[Equation (5.23)]{OB1} showed that the fractional Cauchy problem \eqref{semigroup-0} in dimension $d=1$ with $L_x={\partial^2}/{\partial x^2}$
has solution
 \begin{equation}\label{ob-representation}
m(x,t)=\frac1\gamma \int_0^\infty p(x,u)p_{1/\gamma}(u, \frac1\gamma(2\gamma-1), t)\,du
 \end{equation}
using the parametrization \eqref{char-func}, where $$p(x,u)=T(u)r(x)=\int_{\rr}\frac{e^{-\frac{(x-y)^2}{4u}}}{\sqrt{4\pi u}}r(y)dy$$
 is the heat semigroup corresponding to Brownian motion in $\rr$.  Note that \eqref{ob-representation} involves a stable density, while the equivalent solution \eqref{h1} replaces this by an inverse stable density.  The next result shows how to equate these two forms.  It also extends the result of \cite{OB1} to an abstract fractional Cauchy problem on $\rd$.

\begin{theorem}\label{main-thm}
For $\frac12\leq \gamma <1$, the abstract fractional Cauchy problem \eqref{frac-derivative-0} associated with the Markov process $X(t)$ has two equivalent solutions:
\begin{equation}\begin{split}\label{equiv-representation}
m(x,t)& =\E_x[r(X({E_t}))]=\frac{t}{\gamma}\int_{0}^{\infty}p(x,u)g_{\gamma}(\frac{t}{u^{1/\gamma}})u^{-1/\gamma
-1}\,du \\
&=\E_x[r(X({Y(t)}))|Y(t)>0]=\frac1\gamma \int_0^\infty p(x,u)p_{1/\gamma}(u, \frac1\gamma(2\gamma-1), t)\,du
\end{split}\end{equation}
where $p(x,t)=\E_x[r(X(t))]$ solves the abstract Cauchy problem \eqref{semigroup-0}, $E_t$ is the hitting time \eqref{Edef} of a standard stable subordinator $D(t)$, $g_{\gamma}(x)=p_\gamma(x;\gamma,1)$ is the density of $D(1)$, and $Y(t)$ is a totally negatively skewed stable L\'evy motion with index $\alpha=1/\gamma$ and density $p_{\alpha}(x;2-\alpha, t)$ in the parametrization \eqref{char-func}.
\end{theorem}

\begin{proof}
The integral solution on the first line of \eqref{equiv-representation} was proven in \cite[Theorem 5.1]{limitCTRW}, and then the representation $\E_x[r(X({E_t}))]$ follows from \eqref{Edens}.  This also shows that the integral on the first line of \eqref{equiv-representation} reduces to \eqref{h1}.  Now the second line of \eqref{equiv-representation} follows from Lemma \ref{Ysolve} together with Theorem \ref{YDlink}.
\end{proof}

\begin{remark}
Orsingher and Beghin comment, after equation (5.23) in their paper \cite{OB1}, that the solution to the fractional Cauchy problem \eqref{semigroup-0} in dimension $d=1$ with $L_x={\partial^2}/{\partial x^2}$ and a time derivative of order $\frac12\leq \gamma<1$ can be expressed as  $\E_x[r(B(|Y(t)|))]$ where $B(t)$ is a standard Brownian motion and $Y(t)$ is the stable process from Theorem \ref{main-thm}, so that $u=|Y(t)|$ has the density
$Q(u,t)=\frac1\gamma p_{1/\gamma}(|u|, \frac1\gamma(2\gamma-1), t)$ .   This statement is correct only in the case $\gamma =\frac12$ in which case the process $Y(t)$ is a Brownian motion, and $B(|Y(t)|)$ is an iterated Brownian motion.  In that case, the same result was also recently proven by Baeumer, Meerschaert and Nane \cite{bmn-07}.  However, for $\frac12<\gamma<1$, the spectrally negative stable process $Y$ is not symmetric, and hence the density of $|Y(t)|$ is not equal to $Q(u,t)$, the conditional density of $Y(t)|Y(t)>0$.   In this case, subordination via $|Y(t)|$ does not produce a solution to this fractional Cauchy problem.  Of course, it is still important to find a process whose density equals $Q(u,t)$.  We will return to this question later in the paper, see Remark \ref{TIVO}.
\end{remark}

Meerschaert, Nane, and Vellaisamy \cite{MNV} show that, under some technical conditions, the abstract fractional Cauchy problem \eqref{frac-derivative-0} with $0<\gamma<1$ in a bounded domain $D\subset \rd$ with Dirichlet boundary conditions $m(x,0)=r(x)$ for $x\in D$ and $m(x,t)=0$ for $x\in\partial D$ and $t>0$ is solved by taking
\begin{equation}\label{stoch-rep-L}
m(x,t)=\E_{x}[r(X(E_{t}))I( \tau(X)> E_t)]=\int_{0}^{\infty}p(x,u)h(u,t)du.
\end{equation}
where the first exit time $\tau(X)=\inf \{ t\geq 0:\ X_t\notin D\}$, the uniformly elliptic operator of divergence form $L_x$ is the generator of the semigroup $T(t)f(x)=E_x[f(X_t)I(\tau(X))>t)]$, $p(x,t)=T(t)r(x)$, $E_t$ is the hitting time \eqref{Edef} of the standard stable subordinator, and $h(x,t)$ is the density of $E_t$ as given by \eqref{Edens}.  Then the next result, which extends Theorem \ref{main-thm} to bounded domains, follows immediately from \eqref{EtoY}.

\begin{theorem} Under the technical conditions of \cite[Theorem 3.6]{MNV} the abstract fractional Cauchy problem \eqref{frac-derivative-0} with $1/2\leq\gamma<1$ in a bounded domain $D\subset \rd$ with Dirichlet boundary conditions $m(x,0)=r(x)$ for $x\in D$ and $m(x,t)=0$ for $x\in\partial D$ and $t>0$ has a unique classical solution
\begin{equation}\begin{split}\label{stoch-rep-L-boundary}
m(x,t)&=\E_{x}[r(X(Y(t)))I( \tau(X)> Y(t))|Y(t)>0]\\
&=\alpha\int_0^\infty p(x,u)p_{\alpha}(u, (2-\alpha), t)\,du
\end{split}\end{equation}
with $\tau(X)$, $L_x$, $p(x,t)$ as in the preceding paragraph, $E_t$ is the hitting time \eqref{Edef} of a standard stable subordinator $D(t)$, $g_{\gamma}(x)=p_\gamma(x;\gamma,1)$ is the density of $D(1)$, and $Y(t)$ is a totally negatively skewed stable L\'evy motion with index $\alpha=1/\gamma$ and density $p_{\alpha}(x;2-\alpha, t)$ in the parametrization \eqref{char-func}.
\end{theorem}

\section{Simulation}\label{sectivo}
Recall that the probability density $h(x,t)$ of the inverse stable subordinator $E_t$ defined by \eqref{Edef} solves a time-fractional diffusion equation \eqref{ICPC}.  A simple numerical solution method for \eqref{ICPC} is to simulate a large number of replications of the process $E_t$ and histogram the results.  This method is known as particle tracking \cite{Zhang-PRE2008}, also called the Lagrangian method.  An alternative Eulerian method is based on a finite difference approximation of the fractional derivative \cite{frade}.  In this section, we will examine the implications of space-time duality for both Lagrangian and Eulerian simulation.  Corollary \ref{h-xgov} shows that $h(x,t)$ also solves the space-fractional partial differential equation \eqref{Egovspace} for all $t>0$ and $x>0$.  Generally, space-fractional equations are simpler to simulate.  Lack of memory in time permits efficient Lagrangian simulation based on the underlying Markov process.  The same lack of memory allows Eulerian methods to efficiently step through time.  Since a wide variety of time-fractional partial differential equations can be solved via subordination to the $E_t$ process, the results discussed here have broad applicability.

The space-fractional diffusion equation \eqref{Egovspace} is under-specified, since we desire solutions on $x>0$, and the equivalent equation \eqref{Ygovspace} has another solution $P(x,t)$ supported on the entire real line.  The next result imposes a suitable boundary condition.

\begin{theorem}\label{h-xgov-boundary}
Let $D(1)$ be a stable subordinator with density $g_\gamma(x)=p_\gamma(x;\gamma,b)$ in the parametrization \eqref{char-func}.  Let $E_t$ denote the hitting time process defined by \eqref{Edef}.  Then the density $h(x,t)$ of $E_t$ solves the boundary value problem
\begin{equation}\label{hgovBVP}
\frac{\partial h(x,t)}{\partial t}=b^{-\alpha}\frac{\partial^\alpha h(x,t)}{\partial (-x)^\alpha};\quad\frac{\partial^{\alpha-1}}{\partial(- x)^{\alpha-1} } h(0,t)=0 .
\end{equation}
\end{theorem}

\begin{proof}
Theorem \ref{YDlink} (i) shows that the total mass assigned to the positive real line is $\int_0^\infty P(x,t)\,dx=1/\alpha$, which remains fixed for all $t>0$.  Hence
\begin{equation}\label{fBC}
0=\frac{\partial}{\partial t} \int_0^\infty P(x,t)\,dx=\int_0^\infty \frac{\partial^{\alpha}}{\partial (-x)^{\alpha} } P(x,t)\,dx=\frac{\partial^{\alpha-1}}{\partial(- x)^{\alpha-1} } P(0,t)
\end{equation}
for all $t>0$.  Recall from \eqref{RLnegdef} that fractional derivative ${\partial^{\alpha}}/{\partial(- x)^{\alpha} } P(x,t)$ is defined for $1<\alpha<2$ as the second derivative of the fractional integral of order $2-\alpha$.  Then the last equality in \eqref{fBC} follows from the Fundamental Theorem of Calculus and the fact that ${\partial^{\alpha}}/{\partial(- x)^{\alpha} } P(x,t)$ is defined as the first derivative of that same fractional integral.  Then \eqref{hgovBVP} follows from \eqref{EtoY}.
\end{proof}

\begin{figure}
  \includegraphics[width=4.5in]{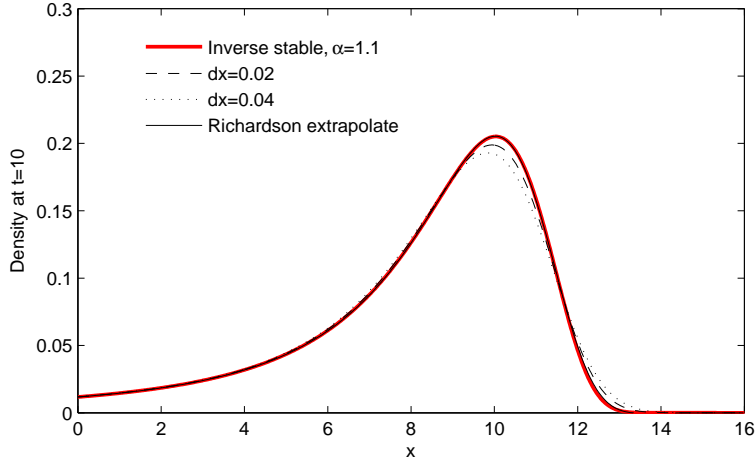}\\
  \caption{Extrapolated Eulerian solution to the boundary value problem \eqref{hgovBVP} matches the inverse stable density $h(x,t)$ from \eqref{Edens}.}\label{FigFD}
\end{figure}

Now the hitting time density $h(x,t)$ can be computed as the point source solution to the space-fractional boundary value problem \eqref{hgovBVP}.  Discretize in space $x_i=i\Delta x$ and time $t_j=j\Delta t$ using the shifted Gr\"unwald finite difference approximation \cite{2sided}:
\[\frac{\partial^\alpha h(x,t)}{\partial (-x)^\alpha}=\lim_{\Delta x\to 0} (\Delta x)^{-\alpha}\sum_{n=0}^{\infty}(-1)^n \binom{\alpha}{n} h(x+(n-1)\Delta x,t),\]
where the fractional Binomial coefficients are defined by
\[\binom{\alpha}{n}=\frac{\Gamma(\alpha+1)}{\Gamma(\alpha-n+1)\Gamma(n+1)} .\]
Then we approximate $h_{ij}=h(x_i,t_j)$ using an explicit Euler scheme
\begin{equation}\label{hEuler}
\frac{h_{ij}-h_{i,j-1}}{\Delta t}=(b\Delta x)^{-\alpha}\sum_{n=0}^{\infty}(-1)^n \binom{\alpha}{n} h_{i+n-1,j-1}
\end{equation}
which reduces to a stable recursive equation for $h_{ij}$, except at the boundary $i=0$ where we apply the boundary condition
\[h_{0j}=-\sum_{n=1}^{\infty}(-1)^n \binom{\alpha-1}{n} h_{i+n-1,j}\]
using the shifted Gr\"unwald approximation once more.  The effect of the boundary condition is to capture the mass that would have exited to the left, into the negative real line, and keep it at $x=0$ to preserve mass.  Figure \ref{FigFD} shows the results of solving boundary value problem \eqref{hgovBVP} in the case $b=1$ via this explicit Euler method.   Figure \ref{FigFD} also shows a semi-analytical solution using \eqref{Edens}, where the stable density $g_\gamma$ was approximated using the algorithm of Nolan \cite{Nolan}.  Richardson extrapolation is based on the fact that the error in the explicit Euler method is approximately proportional to the step size $\Delta x$.  Therefore, a useful estimate of the error at step size $\Delta x$ is $h_{ij}^{2\Delta x}-h_{ij}^{\Delta x}$, and the extrapolated curve is simply the numerical solution at $\Delta x=0.2$ minus this approximate error.

An alternative Lagrangian approach is to simulate the Markov process $t_i=D(x_i)$ at times $x_i=i\Delta x$ as a sum of IID stable random variables with density $g_\gamma(t)=p_\gamma(t;\gamma,b\Delta x)$.  Then we can approximate the inverse process $x_i=E(t_i)$ (at unequally spaced points) and linearly interpolate in $t$.  A histogram of $E(t)$ values from a large number of iterations can be used to approximate the density $h(x,t)$, and more generally, to solve time-fractional diffusion equations via subordination \cite{Zhang-PRE2008}.  An alternative Lagrangian method that requires no interpolation uses the space-time duality from Theorem \ref{YDlink}.  Since $E(t)$ is identically distributed with $Y(t)|Y(t)>0$ we need only simulate the Markov process $Y(t)$ and approximate the conditional density via the histogram.  Note that the proportion of sample paths with $Y(t)>0$ will remain approximately constant since $P(Y(t)>0)=1/\alpha$ for all $t>0$.  Hence this Lagrangian approach is reasonably efficient.

\begin{figure}
  \includegraphics[width=4.5in]{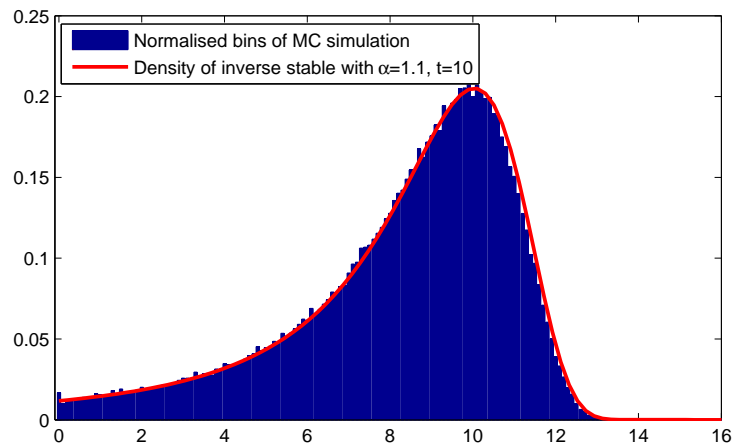}\\
  \caption{Particle tracking simulation of $Z(t)$ matches the conditional density of $Y(t)|Y(t)>0$.}\label{MCTIVO}
\end{figure}

\begin{remark}\label{TIVO}
It is also interesting to find a stochastic process $Z(t)$ whose one dimensional distributions are the same as the conditional distributions of $Y(t)|Y(t)>0$, since the process $Z(t)$ could be used directly as a subordinator to solve time-fractional diffusion equations.  For $\alpha=2$ we can certainly take $Z(t)=|Y(t)|$ by the reflection principle.  This fact is used, for example, to show that iterated Brownian motion $B(Y(t))$ \cite{burdzy1} and Brownian motion in Brownian time $B(|Y(t)|)$ \cite{allouba1,deblassie} have the same one dimensional distributions, and hence the same governing equation.  For $1<\alpha<2$, the process $Y(t)$ is not symmetric, and the reflection principle does not apply.  One possible alternative is to define $Z(t)=Y(u)$ where $u=u(t)$ is the process inverse of
\[t=t(u)=\int_0^u I(Y(s)>0)\,ds\]
so that $t(u)$ is the length of time $Y(s)$ spends being positive during $0<s<u$.  Note that generally $t\leq u$ so that $u=u(t)\geq t$.  In other words, take a sample path of $Y(t)$, snip out the parts where $Y(t)\leq 0$, and glue the remaining parts back together without any gaps in time.  Figure \ref{MCTIVO} compares the results of a particle tracking simulation for this process with the density $h(x,t)$, computed semi-analytically using \eqref{Edens} as in Figure \ref{FigFD}.  In this figure, $b=1$, $\alpha=1.1$, $n=2\times 10^5$ particles were simulated, and a time step of $\Delta t=0.05$ was used in the random walk approximation of $Y(t)$.  The excellent agreement is encouraging.

The process $Z(t)$ is related to local times \cite{bertoin}.  The occupation measure $ \mu_t(B)=\int_0^t I(Y(s)\in B)\,ds$ is well defined for the stable process $Y(t)$, and $\mu_t(dx)$ has a Radon-Nikodym derivative $\ell(x,t)={d\mu_t}/{dx}$ with respect to Lebesgue measure  $dx$ on the real line.  The local time $\ell(x,t)$ measures how much time $Y(s)$ spends at the point $x$ during $0<s<t$.  Now the occupation density formula implies that $t=t(u)=\int_0^\infty \ell(x,u)\,dx$, which can be understood as adding up the time $Y(s)$ spends at all points  $x>0$.  The local time has a scaling property $c^{1 - 1/\alpha}\ell(c^{-1/\alpha}x, \, c^{-1}t)=\ell(x,t)$ in distribution \cite{moduli-cont}, and it follows that $t(cu)=ct(u)$ in distribution.  Then $Z(ct)=c^{1/\alpha}Z(t)$ so that $Z(t)$ has the same scaling as $Y(t)$.  It is known that the inverse local time is a nondecreasing L\'evy process \cite[p.\ 130]{bertoin}.  However, it seems that the integral $t(u)$ is no longer Markovian.  Hence it seems difficult to prove that the pdf of $Z(t)$ is the same as $Y(t)|Y(t)>0$.
\end{remark}

\section{Acknowledgments}

The authors would like to thank Dr. James F. Kelly, Naval Postgraduate School, who originally suggested that Zolotarev duality might be used to connect space and time fractional derivatives.  Thanks also to Professor Yimin Xiao, Department of Statistics and Probability, Michigan State University, for helpful discussions.

\end{document}